\documentclass[numbers=enddot,abstracton]{scrartcl}
\usepackage[utf8]{inputenc}
\usepackage{ccicons}

\usepackage{amsmath}
\usepackage{amscd}

\usepackage{tikz}
\usetikzlibrary{calc}
\usetikzlibrary{arrows}

\DeclareMathAlphabet{\mathpzc}{OT1}{pzc}{m}{it}

\newcommand{\setof}[1]{\{#1\}}
\newcommand{\lam}{\lambda}

\newcommand{\mcP}{\mathcal{P}}

\newcommand{\id}{id}

\newcommand{\colim}[2]{{{\underrightarrow{\lim}}_{#1}{#2}}}

\newcommand{\PP}{\mathbb{P}}
\newcommand{\DD}{D_{(A,B,f)}}
\newcommand{\from}{\leftarrow}
\newcommand{\mcI}{\mathcal{I}}
\newcommand{\mcJ}{\mathcal{J}}
\newcommand{\mcK}{\mathcal{K}}

\newcommand{\upset}[1]{{#1}\!\!\uparrow}
\newcommand{\smupset}[1]{{#1}\!\uparrow}

\newcommand{\eq}{\quad = \quad}

\newtheorem{thm}{Theorem}
\newtheorem{dfn}[thm]{Definition}

\newcommand{\Id}{\mathsf{Id}}

\newcommand{\mcA}{\mathcal{A}}

\renewcommand{\S}{\mathcal S}
\newcommand{\A}{\mathcal{A}}
\renewcommand{\PP}{\mathcal{P}}
\newcommand{\WW}{\mathcal{W}}
\newcommand{\mospace}{\qquad\qquad\!\!\!\!}
\newcommand{\spanme}[5]{\begin{CD} #1 @<#2<< #3 @>#4>> #5 \end{CD}}
\newcommand{\pushout}{\textsf{Pushout}}

\title{On the Local Presentability of $T/\A$}
\author{Patricia Johann and Andrew Polonsky}
\date{July 4, 2018}

\begin{document}

%\Copyright[nc-nd]{Andrew Polonsky and Patricia Johann}

\maketitle

\begin{abstract}
We prove that if $\A$ is a $\lam$-presentable category and $T : \A
\to \A$ is a $\lam$-accessible functor then $T/\A$ is
$\lam$-presentable.
\end{abstract}

\section{Introduction}
Locally presentable categories form a robust class of categories which possess very nice properties, yet are general enough to encompass a large class of examples --- including categories of models of algebraic theories and limit-sketches.

We refer to the standard reference \cite{AR} for definitions and basic properties of locally presentable categories and accessible functors.  These include the following:

\begin{description}
 \item[Proposition] \cite[Prop.1.57]{AR} If $\A$ is a locally $\lam$-presentable category, then for each object $X$, the slice categories $\A / X$ and $X / \A$ are locally $\lam$-presentable.
 \item[Proposition] \cite[Prop.2.43]{AR} If $\A, \A_1$, and $\A_2$ are locally $\lam$-presentable categories and for $i \in \setof{1,2}$, $T_i: \A_i \to \A$ is a $\lambda$-accessible functor, then there exists a regular cardinal $\lam' \ge \lam$ such that the comma category $T_1 \downarrow T_2$ is locally $\lam'$-presentable.
 \item[Proposition] \cite[Exercise 2.h]{AR} If $\A, \A_1, \A_2$ are locally $\lam$-presentable categories and for $i \in \setof{1,2}$, $T_i: \A_i \to \A$ is a $\lambda$-accessible functor \emph{which preserves limits}, then the comma category $T_1 \downarrow T_2$ is locally $\lam$-presentable.
\end{description}

In this note, we add the following to the list above.  (Note that $T/\A$ denotes $T\downarrow\Id$.)
\begin{description}
 \item[Proposition] If $\A$ is a locally $\lam$-presentable category, then for each $\lam$-accessible endofunctor $T : \A \to \A$, the comma categories $\A / T$ and $T / \A$ are locally $\lam$-presentable.
\end{description}

The first claim, that $\A / T$ is $\lam$-presentable, is essentially contained already in the proof of Proposition 2.43 in \cite{AR}.  That proof begins by finding, given $\lam$-accessible $T_1$ and $T_2$, a cardinal $\lam' \ge \lam$ such that $T_1$ and $T_2$ are $\lam'$-accessible \emph{and preserve $\lam'$-presentable objects}.  However, only the fact that $T_1$ preserves such objects is subsequently used.  Since $\Id$ clearly preserves such objects, and $T$ is $\lam$-accessible by hypothesis, we may take $\lam' = \lam$ here and proceed as in \cite{AR} to conclude that $\A / T$ is $\lam$-accessible.

Our main contribution is thus in proving the second claim, that $T / \mcA$ is $\lam$-presentable as well under no additional hypothesis on $T$ beyond $\lam$-accessibility.

\subsection*{Assumptions}
For the rest of this document, we assume the following.

\begin{itemize}
\item $\lambda$ is a regular cardinal.
\item $\A$ is a $\lambda$-presentable category.
\item $\A_0$ is a set of $\lambda$-presentable objects generating all
  of $\A$ under $\lambda$-directed colimits.
\item $T$ is a $\lambda$-accessible endofunctor on $\A$.
\end{itemize}

\section{Outline of the proof}

This section gives a basic summary of our argument. Subsequent
sections will elaborate the individual steps in this
argument. Throughout, if ${\mathcal S} = (S, \leq)$ is a poset
considered as a category, we write $s \in {\mathcal S}$ rather than $s
\in |{\mathcal S}|$. We may also refer to ``the poset $S$'', and
denote the ordering simply by $\leq$, in this case.  If $s \in
{\mathcal S}$ we write $\upset{s}$ for $\{ s' \in {\mathcal S}~|~s
\leq s'\}$.

By the usual arguments, (co)limits in comma categories exist,
including $T/\mcA$, and are computed componentwise. $T/\mcA$ is
therefore cocomplete. The bulk of the argument therefore consists in
exhibiting a set $\PP$ of $\lambda$-presentable objects in $T/\mcA$
that generate all of $T/\mcA$ under $\lambda$-directed colimits.

\begin{description}
\item[1. The set $\PP$.]  Let $\WW = \left\{ (A,P,Q,p,q) \mid A,P,Q
  \in \A_0, p : P \to TA, q : P \to Q \right\}$. For 
\[w = (A,P,Q,p,q) \in \WW\]
let $U(w), f(w)$, and $g(w)$ be defined by the
  pushout
 \begin{equation} \begin{CD}
    P @>p>> TA\\
    @VqVV @VVf(w)V\\
    Q @>g(w)>> U(w)
   \end{CD} \label{e:f}
   \end{equation}

 When we need to refer to individual components of $w \in \WW$, we will write
 $A(w), p(w)$, etc. (For the $w$ above, $A(w)$ is $A$.)

Define 
\[ \PP = \Big\{(A(w),U(w),f(w)) \mid w \in \WW\Big\} \]

Then $\PP$ is a set. Indeed, for every $(A,B,f) \in \PP$ there exists
(at least one) $w \in \WW$ that determines $(A,B,f)$ up to
isomorphism.  We may call such $w$ a \emph{witness} that $(A,B,f) \in
\PP$. Since $\WW$ is clearly a set, there are only set-many witnesses
available.

\item[2. Every element of $\PP$ is $\lambda$-presentable.]  This will
  be proved by a direct argument in Section \ref{s:pres}.
 
\item[3. For every $(A,B,f) \in T/\A$ there is a $\lambda$-directed
  poset $\DD$.] Using the fact that $\A$ is locally $\lambda$-presentable, we first
  collect the following data:
\begin{itemize}
\item Write $A = \colim{i \in \mcI}{A_i}$, where $\mcI$ is a
  $\lam$-directed poset, $\alpha_{i \to i'} : A_i \to A_{i'}$ for $i
  \le i'$, and $\alpha_i : A_i \to A$ for $i \in \mcI$.
\item For each $i \in \mcI$, write $TA_i = \colim{j \in
  \mcJ_i}{P_{i,j}}$, where $\mcJ_i$ is a $\lam$-directed poset,
  $p_{i,j \to j'} : P_{i,j} \to P_{i,j'}$ for $j \le j'$, and $p_{i,j}
  : P_{i,j} \to TA_i$ for $j \in \mcJ_i$.
\item Write $B = \colim{k \in \mcK}{B_k}$, where $\mcK$ is a
  $\lam$-directed poset, $\beta_{k \to k'} : B_k \to B_{k'}$ for $k
  \le k'$, and $\beta_k : B_k \to B$ for $k \in \mcK$.
 \item For each $i \in \mcI$ and each $j \in \mcJ_i$, use the fact
   that $P_{i,j}$ is $\lambda$-presentable, $B = \colim{k \in
     \mcK}{B_k}$, and $f \circ T\alpha_i \circ p_{i,j} : P_{i,j} \to
   B$ to choose $k = k(i,j) \in \mcK$ and $q = q(i,j) : P_{i,j} \to
   B_k$ such that
\[ \begin{CD}
     P_{i,j} @>p_{i,j}>> TA_i\\
     @Vq(i,j)VV @VVf \circ T\alpha_iV\\
     B_{k(i,j)} @>\beta_{k(i,j)}>> B
\end{CD} \]
\end{itemize}

\begin{dfn}
For each $(A,B,f) \in T/\mcA$ define $\DD = (D,\le)$, where
\begin{align*}
 D &\ \ \eq \ \left\{ (i,j,k,q) \quad \Bigg|\quad \begin{aligned}
 &i \in \mcI, j \in \mcJ_i, k \in \mcK, q : P_{i,j} \to B_k\\
 &\qquad f \circ T\alpha_i \circ p_{i,j} = \beta_k \circ q
\end{aligned} \right\} \\[5mm]
\begin{array}{c}
  (i,j,k,q) \\
  \le\\
  (i',j',k',q')
\end{array}  &\;\;\; \Longleftrightarrow \quad
 \left\{\begin{array}{lr}
  i = i, j\le j', k \le k',
  \beta_{k \to k'} \circ q = q' \circ p_{i,j \to j'}\ \  &\fbox{AA}
  \\ 
\hspace*{1.2in}\mbox{or}\\
  i < i', k\le k', \exists r : P_{i,j} \to P_{i',j'}\ \text{such that}\\
  \;\;(T\alpha_{i\to i'} {\circ} p_{i,j} = p_{i',j'} {\circ r}),
   (\beta_{k \to k'} {\circ} q = q' {\circ} r) &\fbox{BB}
  \end{array}
  \right.
\end{align*}
\end{dfn}

As before, we refer to individual components $i,j,k$, and $q$ of $d
{\in} D$ as $i(d), j(d), k(d)$, and $q(d)$. We prove directly that
$\DD$ is a $\lambda$-directed poset in Section~\ref{s:dirpos}.

\item[4. There is a functor $F : \DD \to T/\A$ whose image consists of
  objects in $\PP$\!\!\!].  For $d = (i,j,k,q) \in \DD$, let $w(d) =
  (A_i,P_{i,j},B_k,p_{i,j},q)$.

\begin{dfn}
The functor $F : \DD \to \A$ is defined as follows:
\begin{itemize}
\item $F(d) = (A(w),U(w),f(w))$, where $w = w(d)$ and $A(w)$, $U(w)$,
  and $f(w)$ are as in \eqref{e:f}.
\item $F(d \to d') = (T\alpha_{i \to i'},h) : Fd \to Fd'$ is defined
  whenever $d = (i,j,k,q)$, $d' = (i',j',k',q')$, and $d \leq d'$.
  The morphism $h$ is defined by the pushout property of $U(d)$, as
  shown in either the left or the right diagram below, 
% in Figure~\ref{f:pushout},
  according as $d \le d'$ is obtained via $\fbox{AA}$ or $\fbox{BB}$.
\end{itemize}
\end{dfn}
Note that, for all $d \in \DD$, $F(d)$ is indeed in $\PP$ since $w(d)
\in \WW$.

\begin{figure}[ht]
{\centering
\scriptsize
\begin{tikzpicture}[scale=0.9]
\node (BUL) at (-4, 0) {$P_{i,j}$};
\node (BUR) at (0, 0) {$TA_i$};
\node (BLL) at (-4, -4) {$B_k$};
\node (BLR) at (0, -4) {$U(w(d))$};
\node (FUL) at (-2,-2) {$P_{i,j'}$};
\node (FUR) at (2, -2) {$TA_i$};
\node (FLL) at (-2,-6) {$B_{k'}$};
\node (FLR) at (2,-6) {$U(w(d'))$};
\draw[->,right] (BUL) to node[above]{$\hspace{1.8cm}p_{i,j}$}(BUR);
\draw[->,right] (BLL) to node[above]{\tiny$\hspace{1.8cm}g(w(d))$}(BLR);
\draw[->,right] (FUL) to node[above]{$\hspace{-1.8cm}p_{i,j'}$}(FUR);
\draw[->,right] (FLL) to node[above]{\tiny$\hspace{-1.8cm}g(w(d'))$}(FLR);
\draw[->,right] (BUL) to node[left, yshift=0.8cm]{$q$}(BLL);
\draw[->,right] (BUR) to node[left, yshift=0.8cm]{\tiny$f(w(d))$}(BLR);
\draw[->,right] (FUL) to node[left, yshift=0.8cm]{$q'$}(FLL);
\draw[->,right] (FUR) to node[left, yshift=0.8cm]{\tiny$f(w(d'))$}(FLR);
\draw[->,right] (BUL) to node[above, xshift=0.4cm]{$p_{i,j \to j'}$}(FUL);
\draw[->,right] (BUR) to node[above, xshift=0.4cm]{$\id_{TA_i}$}(FUR);
\draw[->,right] (BLL) to node[above, xshift=0.4cm]{$\beta_{k \to k'}$}(FLL);
\draw[->, dashed] (BLR) to node[above, xshift=0.2cm]{$h$}(FLR);
\node (BUL) at (4, 0) {$P_{i,j}$};
\node (BUR) at (8, 0) {$TA_i$};
\node (BLL) at (4, -4) {$B_k$};
\node (BLR) at (8, -4) {$U(w(d))$};
\node (FUL) at (6,-2) {$P_{i',j'}$};
\node (FUR) at (10, -2) {$TA_i'$};
\node (FLL) at (6,-6) {$B_{k'}$};
\node (FLR) at (10,-6) {$U(w(d'))$};
\draw[->,right] (BUL) to node[above]{$\hspace{1.8cm}p_{i,j}$}(BUR);
\draw[->,right] (BLL) to node[above]{\tiny$\hspace{1.8cm}g(w(d))$}(BLR);
\draw[->,right] (FUL) to node[above]{$\hspace{-1.8cm}p_{i',j'}$}(FUR);
\draw[->,right] (FLL) to node[above]{\tiny$\hspace{-1.8cm}g(w(d'))$}(FLR);
\draw[->,right] (BUL) to node[left, yshift=0.8cm]{$q$}(BLL);
\draw[->,right] (BUR) to node[left, yshift=0.8cm]{\tiny$f(w(d))$}(BLR);
\draw[->,right] (FUL) to node[left, yshift=0.8cm]{$q'$}(FLL);
\draw[->,right] (FUR) to node[left, yshift=0.8cm]{\tiny$f(w(d'))$}(FLR);
\draw[->,right] (BUL) to node[above, xshift=0.4cm]{$r$}(FUL);
\draw[->,right] (BUR) to node[above, xshift=0.4cm]{$T\alpha_{i \to i'}$}(FUR);
\draw[->,right] (BLL) to node[above, xshift=0.4cm]{$\beta_{k \to k'}$}(FLL);
\draw[->, dashed] (BLR) to node[above, xshift=0.2cm]{$h$}(FLR);
\end{tikzpicture}
}
%\caption{The functorial action of $F : \DD \to T/\A$} \label{f:pushout}
\end{figure}

\vspace*{0.1in}

\item[5. There is a cofinal $\lambda$-directed sub(po)set $\DD'$ of
  $\DD$, and thus]
\[\colim{d \in \DD}{Fd} = \colim{d' \in \DD'}{Fd'}\] Let $\DD'$ be the
subposet of $\DD$ comprising $D' = \setof{(i',j',k',q') \in D \mid k'
  \in \upset{k(i',j')}, q' = \beta_{k(i',j') \to k'} \circ q(i',j')}$
under the same ordering as on $D$. To show that $\DD'$ is a cofinal
subposet of $\DD$, suppose $d = (i,j,k,q) \in \DD$.  Note that
 \[ f \circ T\alpha_i \circ p_{i,j} : P_{i,j} \longrightarrow
  \left( \colim{k \in \mcK}{B_k} \right) \] Since $P_{i,j}$ is
  $\lam$-presentable, there exists a $k \in \mcK$ and a morphism $q :
  P_{i,j} \to B_k$ such that $f \circ T\alpha_i \circ p_{i,j} =
  \beta_k \circ q$.  

Since $k(i,j)$ and $q(i,j)$ already satisfy $f \circ T\alpha_i \circ
p_{i,j} = \beta_{k(i,j)} \circ q(i,j)$, by essential uniqueness of
such factorizations there must exist an upper bound $k'$ for $\{k,
k(i,j)\}$ such that
 \begin{align}
    \beta_{k \to k'} \circ q = \beta_{k(i,j) \to k'} \circ q(i,j)
    \label{e:bq}
 \end{align}
Then $d' = (i,j,k',\beta_{k(i,j) \to k'} \circ q(i,j))$ is in $\DD'$,
and \eqref{e:bq} confirms that $d \le d'$ via $\fbox{AA}$.  (Note that
$p_{i,j \to j} = id_{P_{i,j}}$.)

That $\DD'$ is $\lambda$-directed is immediate from it being a
confinal subposet of the $\lambda$-directed set $\DD$.
\clearpage

\item[6. Lemma] For fixed $i$, consider the functor $G_i : \mcJ_i \to
  \A^S$ given by 
\begin{itemize}
\item $G_i(j) = \left(\begin{CD}
      B @<f \circ T\alpha_i \circ p_{i,j}<< P_{i,j} @>p_{i,j}>> TA_i
     \end{CD}\right)$
\item $G_i(j \to j') = \left(\quad
\begin{CD}
  B @<f \circ T\alpha_i \circ p_{i,j}<< P_{i,j} @>p_{i,j}>> TA_i\\
  @V{id_B}VV @Vp_{i,j \to j'}VV @VVid_{TA_i}V\\
  B @<f \circ T\alpha_i \circ p_{i,j'}<< P_{i,j'} @>p_{i,j'}>> TA_i
\end{CD}\quad\right)$
\end{itemize}
Then the following identity holds in $\A^S$:
     \[ \colim{j \in \mcJ_i}{G_i(j)}= \left(
\begin{CD} B @<f \circ T\alpha_i<< TA_i @>id_{TA_i}>> TA_i\end{CD} \right)\]

This is a straightforward computation of a colimit in a functor
category. For completeness we include the proof in Section~\ref{s:lemma}.

\item[7. Lemma] For fixed $i,j$, consider the functor
 $H_{i,j} : \upset{k(i,j)} \to \A^S$ given by
 \begin{itemize}
  \item $H_{i,j}(k) = \left(
  \begin{CD}
    B_k @<\beta_{k(i,j) \to k} \circ q(i,j)<< P_{i,j} @>p_{i,j}>> TA_i
  \end{CD} \right)$.
  \item $H_{i,j}(k \le k') = \left(\qquad
  \begin{CD}
    B_k @<\beta_{k(i,j) \to k} \circ q(i,j)<< P_{i,j} @>p_{i,j}>> TA_i\\
    @V\beta_{k \to k'}VV @Vid_{P_{i,j}}VV @VVid_{TA_i}V\\
    B_k' @<\beta_{k(i,j) \to k'} \circ q(i,j)<< P_{i,j} @>p_{i,j}>> TA_i
  \end{CD}\quad \right)$.
 \end{itemize}
 Then the following identity holds in $\A^S$:
 \[ \colim{k \in \upset{k(i,j)}}{H_{i,j}(k)} = \left(
  \begin{CD} B @<\beta_{k(i,j)} \circ q(i,j)<f \circ TA_i \circ p_{i,j}<
   P_{i,j} @>p_{i,j}>> TA_i
  \end{CD} \right) \]
  
This is also a straightforward computation of a colimit in a functor
category. We omit its proof.

\item[8. Every object $(A,B,f)$ in $T/\A$ is a $\lambda$-directed colimit
  of elements of $\PP$.] We have that $(A,B,f)$ is the colimit of
  $(Fd)_{d \in \DD}$. Indeed, we have

 \begin{align*}
  \colim{d \in \DD}{Fd}
  &\eq \colim{d \in \DD'}{Fd}\\
  &\eq \colim{(i,j,k,q) \in \DD'}{F(i,j,k,q)}\\
  &\eq \colim{i \in \mcI}{\colim{j \in \mcJ_i}{\colim{k \in
        \smupset{k(i,j)}}{F(i,j,k,\beta_{k(i,j) \to k} \circ q(i,j))}}}\\ 
  &\eq \colim{i \in \mcI}{\colim{j \in \mcJ_i}{\colim{k \in \smupset{k(i,j)}}
  {\pushout\left(\spanme{B_k}{\beta_{k(i,j) \to k} \circ
      q(i,j)}{P_{i,j}}{p_{i,j}}{TA_i}\right)}}}\\ 
  &\eq \;\;\;\;\text{\{since pushouts commute with colimits (see Section \ref{s:comm})\}}\\
  & \mospace \colim{i \in \mcI}{\colim{j \in \mcJ_i}
  {\pushout\colim{k \in \smupset{k(i,j)}}
    {\left(\spanme{B_k}{\beta_{k(i,j) \to k} \circ
        q(i,j)}{P_{i,j}}{p_{i,j}}{TA_i}\right)}}}\\ 
  &\eq \;\;\;\;\text{\{by Lemma in point 7\}}\\
& \mospace \colim{i \in \mcI}{\colim{j \in \mcJ_i}{\pushout
  \left(\spanme{B}{f \circ T\alpha_i \circ p_{i,j}}{P_{i,j}}{p_{i,j}}{TA_i}\right)}}\\
  &\eq \;\;\;\;\text{\{since pushouts commute with colimits\}}\\
  & \mospace \colim{i \in \mcI}{\pushout \; \colim{j \in \mcJ_i} 
  {\left(\spanme{B}{f \circ T\alpha_i \circ p_{i,j}}{P_{i,j}}{p_{i,j}}{TA_i}\right)}}\\
  &\eq \;\;\;\;\text{\{by Lemma in point 6\}}\\
  &\mospace \colim{i \in \mcI}{ \pushout
  \left(\spanme{B}{f \circ T\alpha_i}{TA_i}{id_{TA_i}}{TA_i}\right)}\\
  &\eq \;\;\;\;\text{\{pushing out by identity is identity\}}\\
  &\mospace \colim{i \in \mcI}{\left( f \circ T\alpha_i \right)} \\
  &\eq \;\;\;\;\text{\{by computation of colimits in comma categories\}} \\
  & \mospace f
\end{align*}
 This completes the proof of the theorem.
\end{description}

\section{Elements of $\PP$ are $\lam$-presentable}\label{s:pres}
 
Let $D$ be a $\lambda$-directed poset and let $(A^*,B^*,f^*) =
\colim{d \in D}{(A_d,B_d, f_d)}$ in $T/\mcA$. By computation of
colimits in comma categories, we have that
\begin{itemize}
 \item $A^* = \colim{d \in \DD}{A_d}$ in $\mcA$, with structure morphisms
   $\alpha_d : A_d \to A^*$ and $\alpha_{d \to d'}:A_d \to A_{d'}$.
 \item $B^* = \colim{d \in \DD}{B_d}$ in $\mcA$, with structure morphisms
   $\beta_d : B_d \to B^*$ and $\beta_{d \to d'}:B_d \to B_{d'}$.
\item The structure morphisms for $(A^*,B^*,f^*) = \colim{d \in
  \DD}{(A_d,B_d, f_d)}$ in $T/\mcA$ are $(\alpha_d,\beta_d) :
  (A_d,B_d,f_d) \to (A^*,B^*,f^*)$.
 \item Since $T$ is $\lambda$-accessible, $TA^* = \colim{d \in
   \DD}{TA_d}$ in $\mcA$.
\end{itemize}
\noindent
Now let $(A,B,f) \in \mcP$ be determined by $(A,P,Q,p,q)$ via the
pushout
\[ \begin{CD}
    P @>p>> &TA \\
    @VqVV &@VVfV \\
    Q @>g>> &B
   \end{CD}
\] 
and suppose $(\alpha,\beta) : (A,B,f) \to (A^*,B^*,f^*)$ in
$T/\mcA$. Then
\begin{align}
f^* \circ T\alpha = \beta \circ f\label{e:mor1}
\end{align}
\noindent
We want to show that there exists a $d_0 \in \DD$ such that
$(\alpha,\beta)$ factors essentially uniquely through
$(\alpha_{d_0},\beta_{d_0}) : (A_{d_0},B_{d_0},f_{d_0}) \to
(A^*,B^*,f^*)$, as in the characterization of $\lambda$-presentable
objects in Definitions~1.13 and~1.1 of Ad\'{a}mek and Rosick\'{y}
[1994].

Since $\alpha : A \to A^*$, $A$ is $\lambda$-presentable, and $A^*$ is
a colimit, there exists $d \in \DD$, and $\alpha^\circ : A \to A_d$ such that
\begin{equation} \label{e:a}
\alpha = \alpha_d \circ \alpha^\circ
\end{equation}
Similarly, since $\beta \circ g : Q \to B^*$, $Q$ is
$\lambda$-presentable, and $B^*$ is a colimit, there exists $d' \in \DD$
and $g' : Q \to B_{d'}$ such that
\begin{equation}\label{e:g} 
\beta \circ g = \beta_{d'} \circ g'
\end{equation}
Without loss of generality we may assume $d' \ge d$, so that 
\begin{align}
\alpha_d = \alpha_{d'} \circ \alpha_{d \to d'} \label{e:alphas}
\end{align}

Next, observe that
\begin{align*}
\beta_{d'} \circ f_{d'} \circ T\alpha_{d \to d'} \circ T\alpha^\circ \circ p
&= f^* \circ T\alpha_{d'} \circ T\alpha_{d \to d'} \circ T\alpha^\circ
\circ p &&\text{by point 3 above}\\
&= f^* \circ T\alpha_d \circ T\alpha^\circ \circ p &&\text{by
   Equation~\ref{e:alphas}}\\ 
&= f^* \circ T\alpha \circ p &&\text{by Equation~\ref{e:a}}\\
&= \beta \circ f \circ p&&\text{by Equation~\ref{e:mor1}}\\
&= \beta \circ g \circ q&&\text{by definition of $(A,B,f)$}\\
&= \beta_{d'} \circ g' \circ q&&\text{by Equation~\ref{e:g}}
\end{align*}
This exhibits two factorizations of the same morphism from the
$\lambda$-presentable object $P$ to the colimit $B^*$ via
$\beta_{d'}$.  By the essential uniqueness of such factorizations,
there exists a $d_0 \in \DD$, $d_0 \ge d'$ such that
\begin{align} \label{e:y}
\beta_{d' \to d_0} \circ f_{d'} \circ T\alpha_{d \to d'} \circ
T\alpha^\circ \circ p &= \beta_{d' \to d_0} \circ g' \circ q
\end{align}
By the pushout property of $B$, there is therefore a unique
$\beta' : B \to B_{d_0}$ such that
\begin{align}
\beta' \circ f &= \beta_{d' \to d_0} \circ f_{d'} \circ T\alpha_{d \to
   d'} \circ T\alpha^\circ \notag \\ 
&= f_{d_0} \circ T\alpha_{d' \to d_0} \circ T\alpha_{d \to d'} \circ
 T\alpha^\circ \notag \\ 
&= f_{d_0} \circ T\alpha_{d \to d_0} \circ T\alpha^\circ \label{e:bf}\\
\beta' \circ g &= \beta_{d \to d_0} \circ g' \label{e:bg}
\end{align}

Letting $\alpha' = \alpha_{d \to d_0} \circ \alpha^\circ$,
Equation~\eqref{e:bf} states that $(\alpha', \beta') : (A,B,f) \to
(A_{d_0},B_{d_0},f_{d_0})$ is a morphism in $T/\mcA$. That the first
component of this morphism composes with $\alpha_{d_0}$ to $\alpha$ is
obvious:
\[ \alpha_{d_0} \circ \alpha' = \alpha_d \circ \alpha^\circ =
\alpha \] The first equality is by 
Equation~\ref{e:alphas} and the second is by Equation~\ref{e:a}. To
see that the second component of this morphism composes with
$\beta_{d_0}$ to $\beta$, we use uniqueness property of the morphism
from the pushout $B$ to $B^*$. That is, we show that
\begin{align*}
\beta_{d_0} \circ \beta' \circ f = \beta \circ f\\
\beta_{d_0} \circ \beta' \circ g = \beta \circ g
\end{align*}
to conclude that $\beta_{d_0} \circ \beta' = \beta$. To that end,
observe that
\begin{align*}
   \beta_{d_0} \circ \beta' \circ f 
&= \beta_{d_0} \circ f_{d_0} \circ T\alpha_{d \to d_0} \circ T\alpha^\circ &&
   \text{by Equation~\ref{e:bf}}\\
&= f^* \circ T\alpha_{d_0} \circ T\alpha_{d \to d_0} \circ T\alpha^\circ&&
   \text{by point 3 above}\\ 
&= f^* \circ T\alpha_d \circ T\alpha^\circ && \text{by point 1 above}\\ 
&= f^* \circ T\alpha && \text{by Equation~\ref{e:a}}\\
&= \beta \circ f && \text{by Equation~\ref{e:mor1}}\\
   \beta_{d_0} \circ \beta' \circ g 
&= \beta_{d_0} \circ \beta_{d' \to d_0} \circ g' &&
   \text{by Equation~\ref{e:bg}}\\
&= \beta_{d'} \circ g' && \text{by point 2 above}\\
&= \beta \circ g && \text{by Equation~\ref{e:g}}
\end{align*}
So $\beta_{d_0} \circ \beta' = \beta$ and thus $(\alpha,\beta) =
(\alpha_{d_0},\beta_{d_0}) \circ (\alpha',\beta')$ in $T/\mcA$.

To see that this factorization is essentially unique, suppose
$(\alpha, \beta) = (\alpha_{d_0},\beta_{d_0}) \circ
(\alpha'',\beta'')$. We must show that there exists an $l \geq d_0$
such that $(\alpha_{d_0 \to l}, \beta_{d_0 \to l}) \circ (\alpha',
\beta') = (\alpha_{d_0 \to l}, \beta_{d_0 \to l}) \circ (\alpha'',
\beta'')$. Since $(\alpha', \beta')$ and $(\alpha'', \beta'')$ are
morphisms in $T/\mcA$, we have
\begin{align}
\beta' \circ f = f_{d_0} \circ T\alpha' \label{e:com1}\\
\beta'' \circ f = f_{d_0} \circ T\alpha'' \label{e:com2}
\end{align}
Moreover, since $A$ is $\lambda$-presentable and $\alpha_{d_0} \circ
\alpha' = \alpha = \alpha_{d_0} \circ \alpha''$, there exists a $d_1
\geq d_0$ such that
\begin{align}
\alpha_{d_0 \to d_1} \circ \alpha' = \alpha_{d_0 \to d_1} \circ
\alpha'' \label{e:com3} 
\end{align}

Let $\gamma = \beta_{d_0 \to d_1} \circ \beta' \circ q$ and $\gamma' =
\beta_{d_0 \to d_1} \circ \beta'' \circ q$. Then
\begin{align*}
  \beta_{d_1} \circ \gamma 
&= \beta_{d_1} \circ \beta_{d_0 \to d_1} \circ \beta' \circ q && \text{by
    definition of $\gamma$}\\
&= \beta_{d_0} \circ \beta' \circ q\\
&= \beta \circ q && \text{ by hypothesis}\\
&= \beta_{d_0} \circ \beta'' \circ q && \text{ by hypothesis}\\
&= \beta_{d_1} \circ \beta_{d_0 \to d_1} \circ \beta'' \circ q\\ 
&= \beta_{d_1} \circ \gamma' && \text{by
    definition of $\gamma'$}
\end{align*}
That is, we have two factorizations $\beta_{d_1} \circ \gamma$ and
$\beta_{d_1} \circ \gamma'$ of the same morphism from the
$\lambda$-presentable object $Q$ to $A^*$. There must therefore exist
an $l \geq d_1$ such that
\begin{align}
\beta_{d_1 \to l} \circ \gamma = \beta_{d_1 \to l} \circ \gamma' \label{e:com4}
\end{align}

We want to show that 
\[ (\alpha_{d_0 \to l}, \beta_{d_0 \to l}) \circ (\alpha', \beta') =
(\alpha_{d_0 \to l}, \beta_{d_0 \to l}) \circ (\alpha'', \beta'')\]
For the first components we have 
\begin{align*}
  \alpha_{d_0 \to l} \circ \alpha'
&= \alpha_{d_1 \to l} \circ \alpha_{d_0 \to d_1} \circ \alpha'\\
&= \alpha_{d_1 \to l} \circ \alpha_{d_0 \to d_1} \circ \alpha'' \text{ by
    Equation~\ref{e:com3}}\\
&=   \alpha_{d_0 \to l} \circ \alpha''
\end{align*}
For the second components we first observe that
\begin{align*}
 \beta_{d_0 \to l} \circ \beta' \circ f
&= \beta_{d_0 \to l} \circ f_{d_0}\circ T\alpha' && \text{by
   Equation~\ref{e:com1}}\\
&= f_l \circ T\alpha_{d_0 \to l} \circ T\alpha' &&
 \text{$(\alpha_{d_0 \to l},\beta_{d_0 \to l}) :
 (A_{d_0},B_{d_0},f_{d_0}) \to (A_l,B_l,f_l)$}\\
&= f_l \circ T\alpha_{d_0\to l} \circ T\alpha'' && \text{by the
   calculation for the first components}\\
&= \beta_{d_0\to l} \circ f_{d_0}\circ T\alpha'' && 
 \text{$(\alpha_{d_0 \to l},\beta_{d_0\to l}) :
 (A_{d_0},B_{d_0},f_{d_0}) \to (A_l,B_l,f_l)$}\\
&= \beta_{d_0 \to l} \circ \beta'' \circ f && \text{by
   Equation~\ref{e:com2}}\\
\end{align*}
In addition, we have
\begin{align*}
  \beta_{d_0 \to l} \circ \beta' \circ q
&= \beta_{d_1\to l} \circ \beta_{d_0 \to d_1} \circ \beta' \circ q\\
&= \beta_{d_1 \to l} \circ \gamma && \text{ by definition of $\gamma$}\\
&= \beta_{d_1\to l} \circ \gamma' && \text{ by Equation~\ref{e:com4}}\\
&= \beta_{d_1\to l} \circ \beta_{d_0 \to d_1} \circ \beta'' \circ q &&
  \text{ by definition of $\gamma'$}\\ 
&= \beta_{d_0 \to l} \circ \beta'' \circ q
\end{align*}
From these latter two calculations we conclude that $\beta_{d_0 \to l}
\circ \beta' = \beta_{d_0 \to l} \circ \beta''$ by the uniqueness of
the morphism from $B$ to $B^*$ obtained from the fact that $B$ is a
pushout.

 \section{$\DD$ is a $\lam$-directed poset}
 \label{s:dirpos}
 
\begin{description}
\item[Reflexivity.] $(i,j,k,q) \le (i,j,k,q)$ follows
  straightforwardly by $\fbox{AA}$, since $\beta_{k \to k} = id_{B_k}$ and
  $p_{i,j\to j} = id_{P_{i,j}}$, while $\mcI$, $\mcJ_i$, and $\mcK$ are posets. 
 \item[Antisymmetry.] Suppose $(i,j,k,q) \le (i',j',k',q') \le
   (i,j,k,q)$.  Then since $\mcI$ is a poset we have $i \le i' \le i$, so
   that $i=i'$. Similarly, since $\mcK$ is a poset we have $k \le k' \le
   k$ so that $k = k'$.  That $i=i'$ forces both inequalities to arise
   via $\fbox{AA}$, whence the fact that $\mcJ_i$ is a poset entails that
   $j \le j' \le j$ in $\mcJ_i$, so that $j = j'$. Finally, $\fbox{AA}$
   also entails that $q = id_{B_j} \circ q = \beta_{k \to k'} \circ q = q'
   \circ p_{i,j \to j'} = q' \circ id_{P_{i,j}} = q'$, so that
   $(i,j,k,q) = (i',j',k',q')$.
\item[Transitivity.] Suppose $(i,j,k,q) \le (i',j',k',q') \le
  (i'',j'',k'',q'')$. We distinguish four cases.
\begin{description}

\item[\underline{$i=i',i'=i''$}.] Here, $j \le j' \le j''$ yields $j
  \le j''$, and $k \le k' \le k''$ yields $k \le k''$.  Also,
\begin{align*}
  \beta_{k \to k''} \circ q &= \beta_{k' \to k''} \circ \beta_{k \to k'} \circ q\\
  &= \beta_{k' \to k''} \circ q' \circ p_{i,j \to j'}\\
  &= q'' \circ p_{i',j' \to j''} \circ p_{i,j \to j'}\\
  &= q'' \circ p_{i, j \to j''}
\end{align*}
Thus $(i,j,k,q) \le (i'',j'',k'',q'')$ by $\fbox{AA}$.\\

\item[\underline{$i=i',i'<i''$}.] Here, $k \le k' \le k''$ yields $k
  \le k''$ as before, but now 
\begin{align}
   \beta_{k \to k'} \circ q &= q' \circ p_{i,j \to j'}  \label{e:i11}
\end{align}
and there exists an $r' : P_{i',j'} \to P_{i'',j''}$ such that
\begin{align}
  T\alpha_{i'\to i''} \circ p_{i',j'} &= p_{i'',j''} \circ
  r' \label{e:i12}\\ 
  \beta_{k' \to k''} \circ q' &= q'' \circ r' \label{e:i13}
\end{align}
We first observe that $i < i''$. Then we put $r'' = r' \circ p_{i,j \to j'}$
and note that
 \begin{align*}
  T\alpha_{i \to i''} \circ p_{i,j}
  &= T\alpha_{i \to i''} \circ p_{i,j'} \circ p_{i,j \to j'} &&
  \text{ by Definition 1.2}\\
  &= T\alpha_{i' \to i''} \circ p_{i',j'} \circ p_{i,j \to j'} &&
  \text{ since $i = i'$}\\
  &= p_{i'',j''} \circ r' \circ p_{i,j \to j'} && \text{ by Equation~\ref{e:i12}}\\
  &= p_{i'',j''} \circ r'' && \text{ by definition of $r''$}
 \end{align*}
\noindent
and
\begin{align*}
  \beta_{k \to k''} \circ q
  &= \beta_{k' \to k''} \circ \beta_{k \to k'} \circ q && \text{ by
    Definition 1.3}\\
  &= \beta_{k' \to k''} \circ q' \circ p_{i,j \to j'} && \text{ by
    Equation~\ref{e:i11}}\\ 
  &= q'' \circ r' \circ p_{i,j\to j'} && \text{ by
    Equation~\ref{e:i13}}\\ 
  &= q'' \circ r'' && \text{ by definition of $r''$}
 \end{align*}
From this we conclude that $(i,j,k,q) \le (i'',k'',k'',q'')$ by
$\fbox{BB}$.\\

\item[\underline{$i<i',i'=i''$}.]  Here, $k \le k' \le k''$ yields $k
  \le k''$ as before, but now 
\begin{align}
 \beta_{k' \to k''} \circ q' &= q'' \circ p_{i',j' \to j''}  \label{e:i23}
\end{align}
and there exists an $r : P_{i,j} \to P_{i',j'}$ such that
\begin{align}
 T\alpha_{i\to i'} \circ p_{i,j} &= p_{i',j'} \circ r \label{e:i21}\\
 \beta_{k \to k'} \circ q &= q' \circ r \label{e:i22} 
\end{align}
We first observe that $i < i''$. Then we put $r'' = p_{i',j' \to j''}
\circ r$ and note that 
 \begin{align*}
  T\alpha_{i \to i''} \circ p_{i,j}
  &= T\alpha_{i \to i'} \circ p_{i,j} && \text{ since $i' = i''$}\\
  &= p_{i',j'} \circ r && \text{ by Equation~\ref{e:i21}}\\
  &= p_{i',j''} \circ p_{i',j' \to j''} \circ r && \text{ by
    Definition 1.2}\\
  &= p_{i'',j''} \circ r'' && \text{ by definition of $r''$}
\end{align*}
\noindent
and
\begin{align*}
  \beta_{k \to k''} \circ q
  &= \beta_{k' \to k''} \circ \beta_{k \to k'} \circ q && \text{ by
    Definition 1.3}\\
  &= \beta_{k' \to k''} \circ q' \circ r && \text{ by
    Equation~\ref{e:i22}}\\ 
  &= q'' \circ p_{i',j' \to j''} \circ r && \text{ by
    Equation~\ref{e:i23}}\\ 
  &= q'' \circ r'' && \text{ by definition of $r''$}
\end{align*}
From this we conclude that $(i,j,k,q) \le (i'',j'',k'',q'')$ by
$\fbox{BB}$.\\

\item[\underline{$i<i',i'<i''$}.]  Here, $k \le k' \le k''$ yields $k
  \le k''$.  Similarly, $i<i'<i''$ yields $i < i''$.  There exist $r :
  P_{i,j} \to P_{i',j'}$ and $r': P_{i',j'} \to P_{i'',j''}$ such that
\begin{align}
  T\alpha_{i\to i'} \circ p_{i,j} &= p_{i',j'} \circ r \label{e:i31}\\
  \beta_{k \to k'} \circ q &= q' \circ r \label{e:i32} \\
  T\alpha_{i'\to i''} \circ p_{i',j'} &= p_{i'',j''} \circ r' \label{e:i33}\\
  \beta_{k' \to k''} \circ q' &= q'' \circ r' \label{e:i34}
\end{align}
We put $r'' = r' \circ r$ and note that
\begin{align*}
  T\alpha_{i \to i''} \circ p_{i,j}
  &= T\alpha_{i' \to i''} \circ T\alpha_{i \to i'} \circ p_{i,j} && \text{ by
    Definition 1.1}\\
  &= T\alpha_{i' \to i''} \circ p_{i',j'} \circ r && \text{ by 
    Equation~\ref{e:i31}}\\
  &= p_{i'',j''} \circ r' \circ r && \text{ by
    Equation~\ref{e:i33}}\\ 
  &= p_{i'',j''} \circ r'' && \text{ by definition of $r''$}
\end{align*}
\noindent 
and
\begin{align*}
  \beta_{k \to k''} \circ q
  &= \beta_{k' \to k''} \circ \beta_{k \to k'} \circ q && \text{ by
    Definition 1.3}\\
  &= \beta_{k'\to k''} \circ q' \circ r && \text{ by
    Equation~\ref{e:i32}}\\
  &= q'' \circ r' \circ r && \text{ by
    Equation~\ref{e:i34}}\\
  &= q'' \circ r'' && \text{ by definition of $r''$}
\end{align*}
From this we conclude that $(i,j,k,q) \le (i'',j'',k'',q'')$ by
$\fbox{BB}$. 
\end{description}

\item[Directedness.]
Suppose $S \subseteq \DD$, $|S| < \lambda$.  We will construct an upper
bound for $S$ in $\DD$.
 
Let $i^* \in \mcI$ be an upper bound in $\mcI$ for $\setof{i \mid
  (i,j,k,q) \in S}$ and $k_0^* \in \mcK$ be an upper bound in $\mcK$
for $\setof{k \mid (i,j,k,q) \in S}$. For each $s=(i,j,k,q) \in S$, we
will define an element $j(s) \in \mcJ_{i^*}$ and a morphism $r(s) :
P_{i,j} \to P_{i^*,j(s)}$ satisfying
\begin{align}
T\alpha_{i \to i^*} \circ p_{i,j} = p_{i^*,j(s)} \circ r(s) \label{e:rs}
\end{align}
as follows:
\begin{itemize}
\item If $i = i^*$, put $j(s) = j$, and $r(s) = id_{P_{i,j}}$.
\item If $i < i^*$, the facts that $P_{i,j}$ is $\lambda$-presentable
  and
\begin{align*}
T\alpha_{i \to i^*} \circ p_{i,j} 
&: P_{i,j} \to TA_{i^*}\\ 
&= P_{i,j} \to \left(\colim{j' \in J_{i^*}}{P_{i^*,j'}}\right)
\end{align*}
ensure the existence of $j(s) \in \mcJ_{i^*}$ and $r(s) : P_{i,j} \to
P_{i^*,j(s)}$ such that \eqref{e:rs} is valid.
\end{itemize}
Let $j^*$ be an upper bound in $\mcJ_{i^*}$ for the set $\setof{j(s) \mid
  s \in S}$. Now, since $P_{i^*,j^*}$ is $\lambda$-presentable, and
\begin{align*}
f \circ T\alpha_{i^*} \circ p_{i^*,j^*}
&: P_{i^*,j^*} \to B\\
&= P_{i^*,j^*} \to \left(\colim{k \in \mcK}{B_k}\right)
\end{align*}
there also exist $k^*_1 \ge k^*_0$ and $q^*_1 : P_{i^*,j^*} \to B_{k^*_1}$
such that 
\begin{align}
f \circ T\alpha_{i^*} \circ p_{i^*,j^*} = \beta_{k^*_1} \circ
q^*_1 \label{e:q} 
\end{align}

Now, for each $s=(i,j,k,q)$, we have 
\begin{align*}
  \beta_{k^*_1} \circ q^*_1 \circ p_{i^*,j(s) \to j^*} \circ r(s)
  &= f \circ T\alpha_{i^*} \circ p_{i^*,j^*} \circ p_{i^*,j(s) \to
    j^*} \circ r(s) && \text{ by Equation~\ref{e:q}}\\
  &= f \circ T \alpha_{i^*} \circ p_{i^*,j(s)} \circ r(s) && \text{ by
    Definition 1.2}\\
  &= f \circ T \alpha_{i^*} \circ T \alpha_{i \to i^*}
  \circ p_{i,j} && \text{ by Equation~\ref{e:rs}}\\
  &= f \circ T \alpha_i \circ p_{i,j} && \text{ by
    Definition 1.1}\\
  &= \beta_k \circ q && \text{ since $s \in \DD$}\\
  &= \beta_{k^*_1} \circ \beta_{k \to k^*_1} \circ q && \text{ by
    Definition 1.3}
\end{align*}
We therefore have two factorizations of the above morphism from $P_{i,j}$
to $B$ through $\beta_{k^*_1}$. Since $P_{i,j}$ is
$\lambda$-presentable there exists a $k(s) \ge k^*_1$ such that
\begin{align}
  \beta_{k^*_1 \to k(s)} \circ q^*_1 \circ p_{i^*,j(s)\to j^*} \circ r(s)
 = \beta_{k^*_1 \to k(s)} \circ \beta_{k\to k^*_1} \circ q
 = \beta_{k \to k(s)} \circ q \label{eq:q}
\end{align}
Let $k^*$ be an upper bound in $\mcK$ of $\setof{k(s) \mid s \in S}$.

Let $q^* = \beta_{k^*_1 \to k^*} \circ q^*_1$ and define $s^* =
(i^*,j^*,k^*,q^*)$. Then
\begin{align*}
 f \circ T\alpha_i \circ p_{i^*,j^*}
 &= \beta_{k^*_1} \circ q^*_1 && \text{ by Equation~\ref{e:q}}\\
 &= \beta_{k^*} \circ \beta_{k^*_1 \to k^*} \circ q^*_1 && \text{ by
    Definition 1.3}\\
 &= \beta_{k^*} \circ q^* && \text{ by definition of $q^*$}
\end{align*}
so that $s^* \in \DD$.
To see that $s^* \ge s$ for each $s = (i,j,k,q) \in S$, first note
that $i^* \ge i$ and $k^* \ge k^*_0 \ge k$ by construction. Moreover,
\begin{itemize}
\item If $i=i^*$, then $j^* \ge j(s) = j$, and
 \begin{align*}
 \beta_{k \to k^*} \circ q
 &= \beta_{k(s) \to k^*} \circ \beta_{k \to k(s)} \circ q && \text{ by
    Definition 1.3}\\
 &= \beta_{k(s) \to k^*} \circ \beta_{k^*_1 \to k(s)} \circ q^*_1
 \circ p_{i^*,j(s)\to j^*} \circ r(s) && \text{ by
   Equation~\ref{eq:q}}\\  
 &= \beta_{k^*_1 \to k^*} \circ q^*_1 \circ p_{i^*,j(s)\to j^*} \circ r(s)
 && \text{ by Definition 1.3}\\ 
 &= \beta_{k^*_1 \to k^*} \circ q^*_1 \circ p_{i,j \to j^*} \circ
 id_{P_{i,j}} && \text{ since $j = j(s)$ so that $r(s) = id_{P_{i,j}}$}\\
 &= q^* \circ p_{i,j \to j^*} && \text{ by definition of $q^*$}
 \end{align*}
which entails that $(i,j,k,q) \le (i^*,j^*,k^*,q^*)$ by $\fbox{AA}$.
\item If $i<i^*$, put $r = p_{i^*,j(s) \to j^*} \circ r(s) : P_{i,j}
  \to P_{i^*,j^*}$ and observe that
\begin{align*}
T\alpha_{i\to i^*} \circ p_{i,j}
&= p_{i^*,j(s)} \circ r(s) && \text{ by
   Equation~\ref{e:rs}}\\
&= p_{i^*,j^*} \circ p_{i^*,j(s) \to j^*} \circ r(s)  && \text{ by
  Definition 1.2}\\  
&= p_{i^*,j^*} \circ r && \text{ by definition of $r$}
\end{align*}
\noindent
and
\begin{align*}
\beta_{k \to k^*} \circ q
&= \beta_{k(s) \to k^*} \circ \beta_{k \to k(s)} \circ q && \text{ by
  Definition 1.3}\\  
&= \beta_{k(s) \to k^*} \circ \beta_{k^*_1 \to k(s)} \circ q^*_1
\circ p_{i^*,j(s) \to j^*} \circ r(s) && \text{ by
   Equation~\ref{eq:q}}\\
&= \beta_{k^*_1 \to k^*} \circ q^*_1 \circ p_{i^*,j(s) \to j^*} \circ
  r(s) && \text{ by Definition 1.3}\\
&= \beta_{k^*_1 \to k^*} \circ q^*_1 \circ r && \text{ by definition
    of $r$}\\ 
&= q^* \circ r && \text{ by definition of $q^*$}
\end{align*}
Thus $(i,j,k,q) \le (i^*,j^*,k^*,q^*)$ by $\fbox{BB}$.
\end{itemize}
\end{description}

 \section{Proof of Lemma 1}
 \label{s:lemma}
 
 Let $G_i : \mcJ_i \to \A^S$ be as in the statement of Lemma 1.
 
 For concreteness, choose the following presentation for $S$:
 \[ \spanme{0}{v_0}{V}{v_1}{1} \]
 and let $G^*_i = \colim{j \in \mcJ_i}{G_i}$.
 
 Since colimits in functor categories are computed pointwise, we simply verify that
 \begin{itemize}
  \item $G_i^*(0) = \colim{j \in \mcJ_i}{G_i(0)} = \colim{j \in \mcJ_i}{B} = B$
  \item $G_i^*(V) = \colim{j \in \mcJ_i}{G_i(V)} = \colim{j \in \mcJ_i}{P_{i,j}} = TA_i$
  \item $G_i^*(1) = \colim{j \in \mcJ_i}{G_i(1)} = \colim{j \in \mcJ_i}{TA_i} = TA_i$
  \item $G_i^*(v_0) : TA_i \to B$ is the unique map which satisfies,
  for each $j$ in $\mcJ_i$:
  \[ f \circ T\alpha_i \circ p_{i,j} = G^*(v_0) \circ p_{i,j} \]
  Indeed, $f \circ T\alpha_i$ satisfies this constraint, whence $G_i^*(v_0)$ must equal it.
  \item $G_i^*(v_1) : TA_i \to TA_i$ is the unique map which satisfies,
  for each $j$ in $\mcJ_i$:
  \[ p_{i,j} = G_i^*(v_1) \circ p_{i,j} \]
  Indeed, $id_{TA_i}$ satisfies this constraint, whence $G_i^*(v_1)$ must equal it.
 \end{itemize}
 
 \section{Pushouts commute with colimits}
 \label{s:comm}
Recall that a pushout is simply a colimit for a functor whose domain is
the span category
\[ \S \qquad :=\qquad  \cdot \from \cdot \to \cdot \]

Let $F : \DD \to \A^\S$ be a diagram of spans in $\A$.  Then
\begin{align*}
 & \colim{d \in \DD}{\pushout(Fd)}\\
 &\eq \;\;\;\;\text{\{by definition of pushouts\}}\\
 &\mospace \colim{d \in \DD}{\colim{s \in S}{(Fd)s}}\\
 &\eq \;\;\;\;\text{\{by commutativity of colimits\}}\\
 &\mospace \colim{s \in \S}{\colim{d \in \DD}{Fds}}\\
 &\eq \;\;\;\;\text{\{by computation of colimits in functor categories\}}\\
 &\mospace \colim{s \in \S}{\left(\colim{d \in \DD}{Fd}\right)s}\\
 &\eq \pushout\left(\colim{d \in \DD}{Fd}\right)
\end{align*}

\end{document}